\documentstyle[11pt]{article}
\def\C{{\mathbf C}}

\def\mod{{\rm mod}}

\begin{document}

\begin{center}
{\Large \bf Multiple fibers  of holomorphic Lagrangian fibrations}

\bigskip
 {\large \bf Jun-Muk Hwang}\footnote{This work was supported by ASARC.}
{\large \bf and Keiji Oguiso}\footnote{This work was supported by Alexander von Humboldt Foundation and DFG Forschergruppe 790.}

 \end{center}

\bigskip
\begin{abstract} We determine all possible multiplicities of  general singular fibers of a
holomorphic Lagrangian fibration, under the assumption that all
components of the fibers are of Fujiki class. The multiplicities are
at most 6 and the possible values are intricately related to the
Kodaira type of the characteristic cycle.
\end{abstract}



\bigskip

\section{Introduction}

We work in the category of complex analytic spaces. A proper
morphism between two complex analytic spaces is a fibration if it
has connected fibers. Throughout the paper, we will assume that
all components of the fibers of a fibration are of Fujiki class,
i.e., bimeromorphic to compact K\"ahler manifolds.

 Let $(M, \sigma_M)$ be a
holomorphic symplectic manifold and $f: M \to B$ be a fibration
over a complex manifold $B$ whose fibers are Lagrangian with
respect to $\sigma_M$. We say that $f$ is a Lagrangian fibration.
A smooth fiber of $f$ is a complex torus. We do not assume that
$M$ is compact. Our main interest is to understand the structure
of a general singular fiber of $f$. In [HO], a basic structure
theory of a general singular fiber was developed, modulo the
multiplicity of the fiber. A crucial ingredient is the concept of
characteristic cycles, certain connected 1-cycles in the singular
fiber naturally determined by the symplectic form (see Section 2
for the definition). In fact, when the coefficients are divided by
the multiplicity of the fiber, the structure of the characteristic
cycles is exactly the same as the structure of Kodaira's elliptic
singular fibers, with one exceptional case of infinite chain (cf.
Proposition 2.3 and Theorem 2.4). This gives a quite satisfactory
geometric description of a general singular fiber, modulo its
multiplicity. The central remaining question is the possible
values of the multiplicity of a general singular fiber.  This
question was not touched upon in [HO]. The main goal of the
present paper is to give a complete answer to this question.

Here the theory gets somewhat different from Kodaira's theory of
elliptic fibrations, which corresponds to the case of $\dim M =2$.
When $\dim M =2$, except fibers of type $I_m$, a singular fiber
is simply connected and the multiplicity is always 1. So the
theory of multiplicity is quite simple and the detailed geometry
of singular fibers plays little role in the study of the
multiplicity. When $\dim M \geq 4,$ however, a general singular
fiber is not simply connected.  For many types of the
characteristic cycles, the multiplicity can be bigger than 1. The
precise values of the possible multiplicities are intricately
related to the geometry of the singular fiber, as we will see
below. We will determine all the possible values for each type of
the characteristic cycle. Our main result is the following.

\medskip

{\bf Theorem 1.1 (Main Theorem)} {\it Let $f : M \longrightarrow
B$ be a Lagrangian fibration, $D \subset B$ be an irreducible
component of the discriminant divisor, and $b \in D$ be a general
point. Assume that $\dim\, M = 2d \ge 4$. Then
 the multiplicity $n$ of $f^{-1}(b)$ satisfies  $n \le 6$. More precisely,
 let $\Theta$ be the characteristic
cycle with coefficients divided by $n$. Set $\zeta_k = e^{2\pi i/k}$
and denote by $E_{\tau}$ the elliptic curve of period $\tau$. Then
we have the following possibilities.

(i) If $n=6$, then $\Theta$ is isomorphic to the elliptic curve
$E_{\zeta_3}$;

(ii) If $n=5$, then $\Theta$ is  of Kodaira's type $II$;

(iii) If $n=4$, then $\Theta$ is either isomorphic to the elliptic
curve $E_{\zeta_4}$ or of Kodaira's type $IV$;

(iv) If $n=3$, then $\Theta$ is  isomorphic to the elliptic curve
$E_{\zeta_3}$, of Kodaira's type $III$ or of Kodaira's type $I_0^*$;

(v) If $n=2$, then $\Theta$ is  of Kodaira's type $I_{2m}$ ($m \ge
0$), $I_0^*$, $IV$, $IV^*$, or of type $A_{\infty} (=
I_{\infty})$.

In particular, if $\Theta$ is of Kodaira's type $I_{2m-1}$,
$I_m^*$ ($m \ge 1$),
$III^*$ or  $II^*$, then $n =1$, i.e., $f^{-1}(b)$ is not multiple.

 Moreover, all the cases in (i) -(v) are realizable in each dimension
$2d \ge 4$. }

\medskip
In the course of proving our main theorem, we obtain a couple of
results which are of independent interest. In the classification
of characteristic cycles in [HO], a hypothetical case of type
$D_{\infty}$ remained. In Theorem 2.4 below, we remove this case,
thus giving a complete picture of  the classification of
characteristic cycles (modulo multiplicities). Another result of
independent interest is the stable reduction theory presented in
Section 4. This generalizes  Kodaira's stable reduction of
singular fibers of elliptic fibrations to general singular fibers
of  higher dimensional Lagrangian fibrations. In fact, by
exploiting the geometry of the characteristic cycles, the stable
reduction can be constructed from Kodaira's construction for
surfaces.

\medskip
As was the case with [HO], our results have some overlap with
Matsushita's paper [Ma]. He studies the structure of general
singular fibers and their multiplicities under the additional
assumption that the Lagrangian fibration $f:M \to B$ is a
projective morphism.  He used the sophisticated theory of the
toroidal degeneration of abelian varieties, which requires the
assumption that $f$ is projective. Our approach is quite different
and elementary. Other than   the general tools of complex
geometry, all we need is the classical work of Kodaira on
elliptic fibrations. We believe that our approach gives a simpler
and geometrically clearer explanation for both the structure
theory and the classification result, even in the case when $f$ is
projective. Also it should be mentioned that, compared with our
results,  some cases are missing in the classification list of
[Ma].

\section{Review of [HO] with some complements }
\par
\noindent In this section, we recall the definition and basic
properties of characteristic cycles from [HO] and prove some
complementary results.

 Let $M =(M,
\sigma_M)$ be a holomorphic symplectic manifold of dimension $2d
\ge 4$ and $f:M \to B$ be a  Lagrangian fibration over a
$d$-dimensional complex manifold $B$ as defined in the
introduction.  Let $b \in D \subset B$ be a general point of the
discriminant set $D$. By $f^{-1}(b)$ we denote the scheme theoretic fiber
over $b$. We denote by $f^{-1}(b)_{\rm red}$ the reduction of
$f^{-1}(b)$, i.e., the underlying reduced analytic space
of $f^{-1}(b)$. Since we are interested in the local
geometry of $f$ near $b$, we will consider $B$ as a germ of a
neighborhood of $b$. Let us denote by $T$ the complex Lie group
$\C^{d-1}$.

\medskip
{\bf Proposition 2.1} {\it There exists an action of $T$ on $M$
which preserves each fiber of $f$ and the symplectic form
$\sigma_M$, such that the isotropy subgroup at each point of $M$
is discrete.   In particular, the induced action on the dualizing
sheaf $\omega_M$ is trivial, and consequently, the induced action
on the dualizing sheaf of each fiber $f^{-1}(b)$ is trivial. The
singular locus of the reduction $f^{-1}(b)_{\rm red}$ consists
of finitely many $T$-orbits, each of which is a complex torus.}

\medskip
{\it Proof}. This is immediate from  [HO, Proposition 2.2] except
the invariance of the symplectic form. However, the latter is
immediate from the way the $T$-action is constructed. In fact, the
$T$-action is generated by $d-1$ commuting Hamiltonian vector
fields which clearly leaves the symplectic form invariant. $\Box$

\medskip
The following was proved in  [HO, Theorem 1.3].

\medskip
{\bf Proposition 2.2} {\it  For each component $V$ of
$f^{-1}(b)_{\rm red}$, the normalization $\hat{V}$ is a compact
complex manifold of Fujiki class and the Albanese map ${\rm alb}:
\hat{V} \to {\rm Alb}(\hat{V})$ is a fiber bundle whose fiber is
either ${\bf P}_1$ or an ellpitic curve.}

\medskip  Let us call the image of a fiber of the Albanese map in
$f^{-1}(b)_{\rm red}$  a {\it characteristic leaf}.
 We define an equivalence relation on the set
$f^{-1}(b)_{\rm red}$ by declaring two points $x$ and $x'$
equivalent if there is a finite number of characteristic leaves
$C_1, \ldots, C_{\ell}$ such that $x \in C_1$, $x' \in C_{\ell}$ and
$C_{i-1} \cap C_{i} \neq \emptyset$ for all $2 \leq i \leq\ell.$ Let
us call an equivalence class a {\it characteristic curve.} Thus a
characteristic curve $C$ consists of a countable union of
characteristic leaves $\{ C_i \}$. For a characteristic curve $C =
\cup_i C_i$, we define the {\it characteristic cycle} as the
analytic cycle $\Sigma_i  m_i C_i$ where $m_i$ is the multiplicity
of the component of the divisor $f^{-1}(D)$ which contains $C_i$.
The greatest common divisor of $\{m_i\}$ is the multiplicity of the
fiber $f^{-1}(b)$. In [HO, Theorem 1.4] we have shown the following.
We say that the intersection of two smooth curves is
quasi-transversal if the tangent spaces of the curves at the
intersection point are distinct.

\medskip
{\bf Proposition 2.3} {\it The $T$-action on $f^{-1}(b)_{\rm red}$
induces a transitive action on the set of characteristic curves. In
particular, all characteristic curves on $f^{-1}(b)_{\rm red}$ are
isomorphic. Moreover, when $n$ is the multiplicity of the fiber
$f^{-1}(b)$ and $\Sigma_i m_i C_i$ is a characteristic cycle, the
cycle $\Sigma_i \frac{m_i}{n} C_i,$
 is of the form of

(1) one of the singular fibers of a relatively minimal elliptic
fibration listed by Kodaira [Kd, Theorem 6.2];

(2) $1$-cycle of Type $A_{\infty}$, i.e., $1$-cycle $\sum_{i \in
{\mathbf Z}} C_i$ consisting of infinitely many ${\mathbf P}_1$'s
such that $C_i \cap C_{i+1} = \{P_i\}$ (the intersections are
quasi-transversal and $P_i \not= P_j$ if $i \not= j$), and such that
$C_{i} \cap C_{j} = \emptyset$ if $\vert i - j \vert \ge 2$;

(3) $1$-cycle of Type $D_{\infty}$, i.e., $1$-cycle $C_{0} + C_{1} +
\sum_{i \ge 2} 2C_i$ consisting of infinitely many ${\mathbf P}_1$'s
such that $C_{i} \cap C_{i+1} = \{P_i\}$ for each $i \ge 1$, $C_{0}
\cap C_{2} = \{P_0\}$ (all the intersections are quasi-transversal
and $P_i \not= P_j$ if $i \not= j$) and such that $C_{i} \cap C_{j}
= \emptyset$ for other pairs $i \not= j$. }

\medskip
 Examples of characteristic cycles of the type (1) in
Proposition 2.3 are provided by Kodaira's construction. In [HO,
Proposition 4.13], an example belonging to the type (2) was given.
At the time when [HO] was completed, the authors were not aware
whether examples of the type (3) exist or not. It turns out that the
following simple argument excludes the type (3).

\medskip
{\bf Theorem 2.4} {\it Characteristic cycles belonging to the type
(3) in Proposition 2.3 do not exist.}

\medskip
{\it Proof}. Set $F = f^{-1}(b)$ and let $n$ be its multiplicity.
Assume to the contrary that a characteristic cycle $\Theta$ is of
Type $D_{\infty}$. The characteristic cycles are parametrized by
the $T$-action;
$$\Theta_t = n C_{0,t} + n C_{1,t} + \sum_{i \ge 2} 2n C_{i, t}\,\, (t \in T)\,\, .$$
These cycles $\Theta_t$ cover $F$. Consider the Zariski closures
in $F_{\rm red}$:
$$E_1 := \overline{\cup_{t \in T} C_{0, t} \cup C_{1,t}}\, \, ,\,\, E_k := \overline{\cup_{t \in T} C_{k, t}}\,\, (k \ge 2)\,\,
.$$ Then $E_1$ is a finite union of the irreducible components of
$F_{\rm red}$ and each $E_k$ ($k \ge 2$) is an irreducible component of
$F_{\rm red}$. For the multiplicity reason, $E_k \not\subset E_1$
if $k \ge 2$. On the other hand, $F_{\rm red}$ has only finitely
many irreducible components. Thus, there are $\ell > k \ge 2$ such
that $E_{\ell} = E_k$. Choose $k$ that is the minimum among all such
pairs. By definition of the characteristic cycle of Type
$D_{\infty}$, we have $E_{k-1} \cap E_k \neq \emptyset$. Then,
$E_{k-1} \cap E_k$ is the image of the multi-section of the Albanese
map from $\hat{E_k}$ (the normaliazation of $E_k$) and $C_{\ell, t}
\subset E_k$ is the image of a fiber of the Albanese map, under the
normalization map. Thus, $C_{\ell, t}$ meets $E_{k-1} \cap E_k$,
whence, meets $C_{k-1, t'}$ for some $t'$.  This implies that
$E_{\ell-1} = E_{k-1}$. If $k \ge 3$, this contradicts the
minimality of $k$. If $k=2$, then $C_{\ell, t},$ $\ell \ge 3$ meets a
component $C_{1, t'}$ of multiplicity $1$, a contradiction to the
definition of Type $D_{\infty}$.  $\Box$

\medskip
 Now we describe the local structure of the analytic space (scheme)
$f^{-1}(b).$
For an analytic space $F$ and a point $x \in F$, let us denote by
$(F,x)$ the germ of $F$ at $x$. The following is a consequence of
a combination of some results and arguments in [HO].

\medskip
{\bf Proposition 2.5} {\it In the setting of Proposition 2.3, let $x
\in f^{-1}(b)$ be a  point. Then there exist a germ ${\cal R}$ of an
effective divisor in the germ $({\C}^2, 0)$ and an isomorphism of
germs
$$(f^{-1}(b), x) \cong  {\cal R} \times ({\C}^{d-1}, 0),$$  such
that the factor $({\C}^{d-1}, 0)$ is tangent to the $T$-orbits
near $x$. Moreover, the germ of the characteristic cycle $\Theta
\subset f^{-1}(b)$ through $x$ is biholomorphic to the germ of
${\cal R}$ as a cycle in $({\C}^2, 0)$ by the projection $ {\cal
R} \times ({\C}^{d-1}, 0) \to {\cal R}$:
$$
\begin{array}{ccc} (f^{-1}(b), x) &
\stackrel{\cong}{\longrightarrow} & {\cal R} \times ({\C}^{d-1},
0)
\\ \cup & & \downarrow \\ (\Theta, x) &
\stackrel{\cong}{\longrightarrow} & {\cal R}. \end{array}$$ }

\medskip
{\it Proof}. In [HO, Proposition 4.4 (2)], the isomorphism of germs
$$(f^{-1}(b)_{\rm red}, x) \cong {\cal R}_{\rm red} \times ({\C}^{d-1}, 0)$$ is given for
the reduction $f^{-1}(b)_{\rm red}$ with some reduced divisor ${\cal
R}_{\rm red}$ in $(\C^2,0)$. Since $b$ is a general point of the
discriminant set, $(f^{-1}(b),0)$ has the structure of a divisor in
$(\C^{d+1},0)$. Thus if we define the divisor ${\cal R} \subset
(\C^2, 0)$ whose reduction is ${\cal R}_{\rm red}$, by assigning the
multiplicities of its  components to match those of $(f^{-1}(b),0)$,
we get the desired isomorphism
$$(f^{-1}(b), x) \cong  {\cal R} \times ({\C}^{d-1}, 0).$$ Now, to
establish the asserted
  isomorphism of cycles,  it suffices to prove that $(\Theta_{\rm red},x)$
   is biholomorphic to ${\cal R}_{\rm red}$
 by the projection $ {\cal
R}_{\rm red} \times ({\C}^{d-1}, 0) \to {\cal R}_{\rm red}$. From
[HO, Proposition 4.4 (3)], it is clear that  the projection gives
an isomorphism $(\Theta_{\rm red}, x) \cong {\cal R}_{\rm red}$,
except possibly when $\Theta_{\rm red}$ is of Kodaira's type $II$
or $III$. When it is of type $II$, the isomorphism was proved in
the proof of [HO, Proposition 4.7 (2) and (3)]. When it is of type
$III$ and $x$ is the point where two smooth components of
$\Theta_{\rm red}$ intersect tangentially, by [HO, Proposition 4.4
(3)] again, $ (\Theta_{\rm red}, x) \to {\cal R}_{\rm red}$ is
bijective. It suffices to show that the two components of ${\cal
R}_{\rm red}$ has intersection number 2 in $({\bf C}^2, 0)$. This
is immediate from the intersection number consideration in  the
proof of [HO, Proposition 4.11].  $\Box$

\medskip
Proposition 2.5 provides a structure of analytic space (scheme) on a
characteristic cycle. From now on, we will consider characteristic
cycles with this scheme structure.

\section{Multiplicity in the stable case}
\par
\noindent Following [BHPV, Section V.8], we will say that a
characteristic cycle is {\it stable}, if it is   of type $I_b$, $0
\leq b \leq \infty$, modulo the multiplicity,  where $I_{\infty}$
denotes the case (2) in Proposition 2.3. In this section, we will
determine the multiplicity of a general singular fiber in the case
where the characteristic cycle is stable (Proposition 3.5).

\medskip
To start with, we recall some general facts on the multiplicity. Let
$\Delta $ be the unit disk in the complex plane with the origin $0
\in \Delta$ and $h:Z \to \Delta$ be a fibration of a complex
manifold $Z$.  The multiplicity of the fiber $h^{-1}(0)$ is the
largest positive integer which divides the multiplicity (as a
divisor in $Z$) of each component of $h^{-1}(0)$.

\medskip
{\bf Proposition 3.1} {\it Given a fibration of a complex manifold
$h: Z \to \Delta$, suppose the fiber $h^{-1}(b)$ has multiplicity
$n$. Let $\nu: \Delta \to \Delta$ be the cyclic branched cover of
degree $n$ and let $\tilde{Z}$ be the normalization of the fiber
product of $h$ and $\nu$. Then in the natural   commuting diagram
$$\begin{array}{ccc} \tilde{Z} & \stackrel{\pi}{\longrightarrow} & Z \\
\tilde{h} \downarrow & & \downarrow h \\ \tilde{\Delta} &
\stackrel{\nu}{\longrightarrow}  &\Delta, \end{array} $$ $\pi$ is
an unramified covering of degree $n$ and the fiber
$\tilde{h}^{-1}(0)$ has multiplicity 1.  The  cyclic Galois group
of order $n$ acts on the fiber $\tilde{h}^{-1}(0)$ freely. Assume
further that the dualizing sheaf $\omega_Z$ is trivial, $\omega_Z
\cong {\cal O}_Z$. Then the induced action of the cyclic Galois
group on the dualizing sheaf $\omega_{\tilde{Z}}$ is trivial,
i.e., the action has weight  1.  }

\medskip
{\it Proof}. This is essentially [BHPV, Chapter III, Proposition
9.1] where it is stated and proved when $\dim Z =2$. The same
proof works in any dimension. $\Box$

\medskip
Recall that we have provided a characteristic cycle $\Theta$ with
the structure of a complex analytic space (cf. the remark after
Proposition 2.5). Denote by ${\rm Aut}(\Theta)$ the biholomorphic
automorphism group of $\Theta$ and by ${\rm Aut}^{\omega}(\Theta)$
the subgroup consisting of the automorphisms acting trivially on
the dualizing sheaf $\omega_{\Theta}$. Then the quotient group
${\rm Aut}(\Theta)/{\rm Aut}^{\omega}(\Theta)$ acts faithfully on
$\omega_{\Theta}$.

\medskip
{\bf Proposition 3.2} {\it  Let $F$ be a  general singular fiber
of a Lagrangian fibration with multiplicity 1.  We are given a
$T$-action on $F$ by Proposition 2.1. Let ${\rm Aut}_T(F)$ be the
group  of automorphisms of $F$ commuting with the $T$-action. Let
${\rm Aut}_T^{\omega}$ be the subgroup of ${\rm Aut}_T(F)$
consisting of automorphisms acting trivially on $\omega_F$. Fix a
characteristic cycle $\Theta \subset F$. Then there exists a
canonical injective homomorphism
$$ \theta: {\rm Aut}_T(F)/{\rm Aut}^{\omega}_T(F) \to {\rm
Aut}(\Theta)/{\rm Aut}^{\omega}(\Theta)$$ such that the eigenvalue
of an element $g \in {\rm Aut}_T(F)$ on $\omega_F$ agrees with that
of $\theta(g)$ on $\omega_{\Theta}$.}

\medskip
{\it Proof}.
  Given $g \in {\rm
Aut}_T(F)$ and a fixed choice $\Theta \subset F$, there exists $t
\in T$ such that $t \cdot \Theta = g \cdot \Theta$. Then $t^{-1}
\circ g: \Theta \to \Theta$ is an automorphism of $\Theta$. Define
$$\theta(g) := t^{-1} \circ g \;\; {\rm mod } \; {\rm Aut}^{\omega}(\Theta)
\; \in \; {\rm Aut}(\Theta)/{\rm Aut}^{\omega}(\Theta).$$ We claim
that this definition is independent of the choice of $t$. In fact,
if $(t_1 \circ t_2) \Theta = \Theta$, then $t_1 \circ t_2$
regraded as an automorphism of $\Theta$ is in ${\rm
Aut}^{\omega}(\Theta)$. This is because $T$-action on $F$
gives a trivial action on $\omega_F$ by Proposition 2.1. Since ${\rm
Aut}_T(F)$ commutes with the $T$-action, it is immediate that
$\theta$ is a group homomorphism. Let us determine the kernel of
$\theta$. Suppose that $g$ is in the kernel of $\theta$. Then
$t^{-1} \circ g \in {\rm Aut}^{\omega}(\Theta).$ Since $T$ action
on $\omega_F$ is trivial, this means that $g$ acts trivially on
$\omega_F$, i.e., $g \in {\rm Aut}_T^{\omega}(F)$. It follows that
$\theta$ descends to an injection
$$ \theta: {\rm Aut}_T(F)/{\rm Aut}^{\omega}_T(F) \to {\rm
Aut}(\Theta)/{\rm Aut}^{\omega}(\Theta).$$ From the definition of
$\theta$, it is clear that the eigenvalue of an element $g \in
{\rm Aut}_T(F)$ on $\omega_F$ agrees with that of $\theta(g)$ on
$\omega_{\Theta}$. $\Box$

\medskip
Now let $\Theta$ be a characteristic cycles of type $I_b, 0 \leq b
\leq \infty.$ More precisely, we have the following description.

(1) When $b=0$, $\Theta$ is an elliptic curve.

(2) When $1 \leq b < \infty$, $\Theta$ is the Kodaira fiber of
type $I_b$, consisting of $b$ smooth rational curves.

(3) When $b = \infty$, $\Theta$ is the analytic space with
infinitely many irreducible components described in (2) of
Proposition 2.3.

\medskip
{\bf Proposition 3.3} {\it Let $\Theta$ be of type $I_b, 0 \leq
b\leq \infty$.  Then ${\rm Aut}(\Theta)/{\rm
Aut}^{\omega}(\Theta)$ is a finite cyclic group. Its order is  6
when $\Theta \cong E_{\zeta_3}$, 4 when $\Theta \cong
E_{\zeta_4}$, and  2 otherwise. When $b$ is odd, an element $\tau
\in {\rm Aut}(\Theta)$ which is not in $ {\rm
Aut}^{\omega}(\Theta)$ fixes a singular point $P$ of $\Theta$ and
exchanges the two irreducible components of the germ of $\Theta$
at $P$. }

\medskip
{\it Proof}. The statement  for $I_0$ and $I_1$ is well-known fact
from the theory of elliptic curves (see for instance [Mc,
Proposition 4.2] for the explicit statement with a proof). The
statement for  $I_b, b>0$ then follows from the fact that $\Theta$
of type $I_b$ is an unramified cyclic Gorenstein covering of $I_1$.
 $\Box$

\medskip
{\bf Proposition 3.4} {\it Let $F$ be a general singular fiber of a
Lagrangian fibration with multiplicity 1 whose characteristic cycle
is  of type $I_b$. Suppose that $\Gamma \subset {\rm Aut}_T(F)$ acts
faithfully on $\omega_F$. Then the order of $\Gamma$ is  $1$ or $2$
when $ b >0$ and  $ 2,3, 4$ or $ 6 $ when $ b =0$. When it is 3 or
6, the characteristic cycle is isomorphic to $E_{\zeta_3}$ and when
it is 4, the characteristic cycle is isomorphic to $E_{\zeta_4}$.}

\medskip
{\it Proof}.  By the assumption, $\Gamma \to {\rm Aut}_T(F) / {\rm
Aut}_T^{\omega}(F)$ is injective. Thus the result follows from
Proposition 3.2 and Proposition 3.3. $\Box$

\medskip
{\bf Proposition 3.5} {\it The multiplicity of a general singular
fiber of a Lagrangian fibration with stable characteristic cycles of
 type $I_b$ (modulo multiplicity) can take only the following values:
$1$ or $2$ when $b \geq 1$ and $b$ is even, $1$ when $b$ is odd,
and $1,2,3,4,$ or $6$ when $b =0$. When it is 3 or 6, the reduced
characteristic cycle is isomorphic to $E_{\zeta_3}$ and when it is
4, the reduced characteristic cycle is isomorphic to $E_{\zeta_4}$.}

\medskip
{\it Proof}. Let $f:M \to B$ be a holomorphic Lagrangian fibration
over a complex manifold $B$.  Let $D \subset B$ be a component of
the discriminant of $f$ and $b \in D$ be a general point. Choose an
arc $\alpha: \Delta \to B$ with $\alpha(0) = b $ which intersects
$D$ transversally and denote by $h:Z \to \Delta$
 the pull-back of $f$ by $\alpha$. Then $Z$ is a complex manifold
 and the multiplicity of the fiber $h^{-1}(0)$ is the same as the
 multiplicity of $f^{-1}(b)$. The $T$-action on $M$ induces a $T$-action on $Z$ preserving
 the fibers of $h$.  By applying Proposition 3.1, we get an
 unramified cover $\tilde{Z}$ of $Z$ with a fibration $\tilde{h}:
 \tilde{Z} \to \Delta$. The $T$-action on $Z$ lifts to a $T$-action on $\tilde{Z}$
 and the action of the cyclic Galois group on $\tilde{Z}$ commutes with this
 $T$-action. Since $\omega_Z$ is trivial, the Galois group acts on
 $\omega_{\tilde{Z}}$ trivially. But the action on $\omega_{\Delta}$
 is faithful. Thus the Galois group acts faithfully on the dualizing
sheaf of the fiber $\tilde{h}^{-1}(0)$.
 From the way it is constructed, the fiber $F := \tilde{h}^{-1}(0)$ can be
 realized as a general singular fiber of multiplicity 1 in some
 holomorphic Lagrangian fibration, say, $\tilde{M} \to \tilde{B}$,
where $\tilde{B}$ is the cyclic covering of order
$n$ branched along $D$ and $\tilde{M}$ is the normalization
of the fiber product, near over $b$. Thus we can apply Proposition 3.4
 with the cyclic Galois group as the group $\Gamma$.
This proves Proposition 3.5 except that the order of $\Gamma$ is
not $2$ when $b$ is odd.

To complete the proof, assuming that $b$ is odd and $\Gamma =
\langle \iota \rangle$ is of order $2$, we shall derive a
contradiction. Put $F = \tilde{h}^{-1}(0)$ and choose a
characteristic cycle $\Theta \subset F$. $\Theta$ is either of
type $I_b$ or of Type $I_{2b}$.

If it is of type $I_{2b}$, then $\iota(\Theta) = \Theta$. The
action is a cyclic permutation of the components of $\Theta$. This
is because the characteristic cycle of $h$ is also of type $I_b$.
However, then $\iota^*\omega_{\Theta} = \omega_{\Theta}$, a
contradiction to Proposition 3.2.

Consider the case where $\Theta$ is of type $I_b$. By Proposition
3.2 and its proof, there is $t \in T$ such that $t^{-1} \circ
\iota(\Theta) = \Theta$. Put $\tau = t^{-1} \circ \iota$. Then
$\tau^*(\omega_F) = -\omega_F$. This is because $\iota$ acts on
the base as $-1$ but the action on $\omega_{\tilde{Z}}$ is
trivial. Then, by Proposition 3.2, $\tau^*\omega_{\Theta} =
-\omega_{\Theta}.$ As $b$ is odd, by Proposition 3.3,  $\tau$
fixes a singular point $P$ of $\Theta$ and exchanges the two
components, say $C_1$ and $C_2$, of the germ of $\Theta$ at $P$.
 Now consider the isotropy action of $\tau$ on the
$(d+1)$-dimensional tangent space ${\rm T}_P(\tilde{Z})$  of
$\tilde{Z}$ at $P$. Because $\iota$ commutes with $T$, the
isotropy action fixes the $(d-1)$-dimensional subspace ${\rm
T}_{P}(T\cdot P) \subset {\rm T}_{P}(\tilde{Z})$, the tangent to
the orbit $T \cdot P$. Let $N = {\rm T}_{P}(\tilde{Z})/{\rm
T}_{P}(T\cdot P)$ be the two-dimensional quotient space. The
induced action on $N$ is an order-2 automorphism of $N$ which
interchanges two distinct subspaces of dimension 1 in $N$
corresponding to the tangent directions of $C_1$ and $C_2$ at $P$.
Any automorphism of order 2 on a 2-dimensional vector space with
this property must act on $\wedge^2 N$ by $-1$. This implies
$\tau^* \omega_{\tilde{Z}} = -\omega_{\tilde{Z}}$, a
contradiction.
 $\Box$

\section{Stable reduction}
\par
\noindent  When $\Theta$ is unstable, the group ${\rm
Aut}(\Theta)/{\rm Aut}^{\omega}(\Theta)$ is no longer as simple as
in Proposition 3.3. To bound the multiplicity in this case, we need
the stable reduction of the unstable fiber. The goal of this section
is to explain the construction of the stable reduction.

\medskip
We start with recalling Kodaira's stable reduction for a minimal
elliptic fibration  with unstable singular fiber ([Kd, Sections
8,9]). Let $j: S \to \Delta$ be a minimal elliptic fibration with
a singular fiber $j^{-1}(0)$ of unstable type.  There are three
different cases: (Case 1) Type $I^*_b, b \geq 1$, (Case 2) Type
$I^*_0, II^*, III^*, IV^*$ and (Case 3) Type $II, III, IV$. Let us
recall the construction of the stable reduction in each case. A
good reference is [BHPV, Section V.10].

\medskip
(Case 1) Type $I^*_b, b \geq 1$.   There are four $(-2)$-curves in
the central fiber $j^{-1}(0)$. Contracting these four curves gives
a normal surface $S'$ with a fibration $j': S' \to \Delta$ whose
central fiber is a string of $b+1$ nonsingular rational curves of
multiplicity $m=2$. Take a cyclic cover of order $m=2,$ $\nu:
\Delta \to \Delta$ and let $j^{\sharp}: S^{\sharp} \to \Delta$ be
the normalization of the fiber product of $j$ and $\nu$. Then
$S^{\sharp}$ is non-singular and the central fiber of $j^{\sharp}$
is a reduced fiber of type $I_{2b}$. This $j^{\sharp}$ is the
stable reduction of $j$. The Galois group $\langle -1 \rangle$
of order $2$ acts on
$S^{\sharp}$ with 4 fixed points. The generator $-1$ of the Galois
group acts on $\omega_{j^{\sharp-1}(0)}$ by $-1$. The quotient surface
$S^{\sharp}/\langle -1 \rangle$ has 4 singular points of type
$\frac{1}{2}(1,1)$ (under the notation [BPHV, Section V, 10]). From the
construction, we can choose a coordinate system $(x, u)$ at a
fixed point, say $P$, and a coordinate system $(y, v)$ on an open subset in
$S$ where $x$ (resp. $y$) is the pull-back of the coordinate on
$\Delta$ by $j^{\sharp}$ (resp. $j$) such that there exists a
meromorphic correspondence
$$y = x^2, \; v = \frac{ux}{x^2} = \frac{u}{x}.$$
Here $v$ is an affine coordinate of the exceptional curve arising from the
minimal resolution of $S^{\sharp}/\langle -1 \rangle$ at the image
of $P$. The dualizing sheaves of
$S^{\sharp}$ and $S$ are related by
$$dy \wedge dv = 2 dx \wedge du.$$

\medskip
(Case 2) Type $I^*_0, II^*, III^*, IV^*$  There is one irreducible
component $C$ of maximal multiplicity $m$ in $j^{-1}(0)$. The
values of $m$ are (cf. Line 4 of Table 5 in [BHPV, Section V.10]):
$$ m=2 \mbox{ for Type } I^*_0, \; m=6 \mbox{ for Type } II^*, \;
m=4 \mbox{ for Type } III^*, \; m=3 \mbox{ for Type } IV^*.$$ The
connected components of $j^{-1}(0) \setminus mC$ are four
$(-2)$-curves for Type $I^*_0$ and three Hirzebruch-Jung strings
for Type $II^*, III^*, IV^*$. Contract these Hirzebruch-Jung
strings to get a fibration $j_2: S_2 \to \Delta$ of a normal
surface $S_2$ with irreducible central fiber of multiplicity $m$.
Then take a cyclic cover $\nu: \Delta \to \Delta$ of order $m$ and
let $j^{\sharp}: S^{\sharp} \to \Delta$ be the normalization of
the fiber product of $\nu$ and $j_2$. Then $j^{\sharp}$ is a
non-singular elliptic fibration. This $j^{\sharp}$ is the stable
reduction of $j$. The cyclic Galois group of order $m$ acts on
$S^{\sharp}$ with fixed points. The generator $\zeta_m$ of the
Galois group acts on the dualizing sheaf
$\omega_{j^{\sharp-1}(0)}$ by $\zeta_m^{m-1}$ (cf. Line 7 of Table
5 in [BHPV, Section V.10]) while it acts on $\omega_{\Delta}$ by
$\zeta_m$. Consequently, it acts on $\omega_{S^{\sharp}}$
trivially.
 From the construction, $\zeta_m$ has a fixed point, say $P$,
and the quotient surface $S^{\sharp}/\langle \zeta_m \rangle$
has a singular point of type $\frac{1}{m}(1,-1)$ at the image
of $P$. So, we can
choose a coordinate system $(x, u)$ at a fixed point of
$S^{\sharp}$ and a coordinate system $(y, v)$ on an open subset in
$S$ where $x$ (resp. $y$) is the pull-back of the coordinate on
$\Delta$ by $j^{\sharp}$ (resp. $j$) such that there exists a
meromorphic correspondence
 $$y = x^m, \; v = \frac{ux}{x^m} = \frac{u}{x^{m-1}}.$$
Here $v$ is again an affine coordinate of one of the exceptional
curves arising from the minimal resolution $S$ of
$S^{\sharp}/\langle \zeta_m \rangle$ at the image of $P$. The
dualizing sheaves of $S^{\sharp}$ and $S$ are related by
$$dy \wedge dv = m  dx \wedge du.$$

(Case 3) Type $II, III, IV$.  First apply a finite number of
blow-ups of $S$ to get $j_1: S_1 \to \Delta_1$ such that the
central fiber has normal crossing support. The blow-ups needed in
this process are listed in Page 209 of [BHPV, Section V.10]. In
each step, the blow-up center is the unique singular point of the
reduced fiber.  There is one component of maximal multiplicity $m$
in $j_1^{-1}(0)$. The values of $m$ are (cf. Line 4 of Table 5 in
[BHPV, Section V.10]): $$
 m=6 \mbox{ for Type } II, \; m=4 \mbox{ for Type } III, \; m=3 \mbox{
for Type } IV.$$ Beside this component with multiplicity $m$,
there are three smooth rational curves with negative
self-intersection in the singular fiber of $j_1^{-1}(0)$. So they
are Hirzebruch-Jung strings of length 1. Contract these
Hirzebruch-Jung strings to get a fibration $j_2: S_2 \to \Delta$
of a normal surface $S_2$ with irreducible central fiber of
multiplicity $m$. Then take a cyclic cover $\nu: \Delta \to
\Delta$ of order $m$ and  let $j^{\sharp}: S^{\sharp} \to \Delta$
be the normalization of the fiber product of $\nu$ and $j_2$. Then
$j^{\sharp}$ is a non-singular elliptic fibration. This
$j^{\sharp}$ is the stable reduction of $j$. The cyclic Galois
group $\langle \zeta_m \rangle$ of order $m$ acts on $S^{\sharp}$
with a fixed point, say $P$. The generator $\zeta_m$ of the Galois
group acts on the dualizing sheaf $\omega_{j^{\sharp-1}(0)}$ by
$\zeta_m$ (cf. Line 7 in Table 5 in [BHPV, Section V.10]) while it
acts on $\omega_{\Delta}$ by $\zeta_m$. Consequently, it acts on
$\omega_{S^{\sharp}}$ by $\zeta_m^2$ and the quotient surface
$\omega_{S^{\sharp}}/\langle \zeta_m \rangle$ has a singular point
of type $\frac{1}{m}(1,1)$ at the image of $P$. From this
description, we can choose a coordinate system $(x, u)$ at a fixed
point $P$ of $\zeta_m$ and a coordinate system $(y, v)$ on an open
subset in $S$ where $x$ (resp. $y$) is the pull-back of the
coordinate on $\Delta$ by $j^{\sharp}$ (resp. $j$) such that there
exists a meromorphic correspondence
 $$y = x^m, \; v = \frac{ux^{m-1}}{x^m} = \frac{u}{x}.$$
Again, as before, $v$ is an affine coordinate of the exceptional
curve of the minimal resolution $S_1$ of
$\omega_{S^{\sharp}}/\langle \zeta_m \rangle$ at the image of $P$.
Note that this exceptional curve certainly survives under taking
the relative minimal model over the base, of the minimal
resolution of $\omega_{S^{\sharp}}/\langle \zeta_m \rangle$. Thus
$(y,v)$ can be regarded as a coordinate system on an open subset
in $S$. The dualizing sheaves of $S^{\sharp}$ and $S$ are related
by
$$dy \wedge dv = m x^{m-2} dx \wedge du.$$

\medskip
To generalize this construction  to higher dimensions, it is
convenient to introduce the following notion.

\medskip
{\bf Definition 4.1} {\it Let $j: S \to \Delta$ be a fibration of
a normal surface $S$ over a disk which is smooth over $\Delta
\setminus \{0\}.$  A $(d+1)$-dimensional normal complex analytic
variety $Z$ and a fibration $h:Z \to \Delta$ which is smooth over
$\Delta \setminus \{0\}$,  is called a fibration {\it modeled on}
$j: S \to \Delta$ if the following holds.

(i) The Lie group $T= \C^{d-1}$ acts on $Z$  by a holomorphic map
$\gamma: Z \times T \to Z$  which preserves the fibers of $h$,
i.e., the following diagram commutes, where $\bar{h}$ is the
composition of $h$ and the projection to $Z$.
$$
\begin{array}{ccc} Z \times T & \stackrel{\gamma}{\longrightarrow}
& Z \\ \bar{h} \downarrow & & \downarrow h \\ \Delta &
\stackrel{=}{\longrightarrow} & \Delta,
\end{array}.$$ The stabilizer of the $T$-action at each point of $Z$ is discrete and
the singular loci of $h^{-1}(0)_{\rm red}$ consists of finitely
many $T$-orbits, which are $(d-1)$-dimensional tori.

(ii) For each point $x \in h^{-1}(0)$, there exists a
1-dimensional compact connected analytic subscheme $\Psi_x \subset
h^{-1}(0)$ containing $x$ such that for $x, y \in h^{-1}(0)$,
either $\Psi_x \cap \Psi_y = \emptyset$ or $\Psi_x = \Psi_y$, and
for $g \in T$,
$$ g \cdot \Psi_x = \Psi_{g \cdot x}, \mbox{ i.e., }
\gamma( \Psi_x, g) = \Psi_{\gamma(x,g)}.$$

(iii) For each $x \in h^{-1}(0)$, there exist a point $s \in
j^{-1}(0) \subset S$ and a morphism $$\iota: j^{-1}(0) \to Z, \;
\iota(s) =x,$$ inducing an isomorphism $j^{-1}(0) \cong \Psi_x.$
Furthermore, there exists  isomorphisms of germs $\rho: (S,s)
\times (T, 0) \cong (Z,x)$ and $\rho|_{j^{-1}(0) \times T}:
(j^{-1}(0), s) \times (T,0) \cong (h^{-1}(0), x)$ with the
commuting diagrams
$$\begin{array}{ccccccc}   (S,s) \times (T,0) & \stackrel{\rho}{\longrightarrow} & (Z,x)
& & (j^{-1}(0),s) \times (T,0) &
\stackrel{\rho|_{j^{-1}(0)\times T}}{\longrightarrow} & (h^{-1}(0),x)   \\
j \times \{0\}
 \downarrow & & \downarrow h & &  \cap & & \cap   \\
\Delta & \stackrel{=}{\longrightarrow} & \Delta & & (S,s) \times
(T,0) & \stackrel{\rho}{\longrightarrow} & (Z,x),\end{array}$$
compatible with the $T$-action $$\begin{array}{ccc} (j^{-1}(0),s)
\times (T,0) & \stackrel{\iota \times {\rm
Id}_T}{\longrightarrow} & (Z,x) \times (T,0)  \\  \rho \downarrow & & \downarrow \gamma \\
(Z,x) & \stackrel{=}{\longrightarrow} & (Z,x).\end{array}$$ }

\medskip
{\bf Proposition 4.2} {\it In the setting of Definition 4.1, assume
that $h^{-1}(0)_{\rm red}$ is not smooth. Then there exists a
complex torus $T'$ of dimension $d-1$ with a $T$-action and a
$T$-equivariant morphism $q: h^{-1}(0) \to T'$ such that for each $t
\in T'$, $q^{-1}(t)$ is $\Psi_x$ for some $x$.}

\medskip
{\it Proof}. For a given $\Psi_x, x \in h^{-1}(0)$, let
$I_{\Psi_x}$ be the ideal sheaf on $h^{-1}(0)$ defining $\Psi_x$
as a subscheme. By Definition 4.1 (iii), a germ of $h^{-1}(0)$ is
 the product of a germ of $\Psi_x$ with a smooth germ of dimension $(d-1) = \dim T$.
 Thus  the conormal sheaf
$I_{\Psi_x}/I^2_{\Psi_x}$ is a locally free sheaf on $\Psi_x$ of
rank $d-1$. Moreover, the vector fields generating $T$-action on
$h^{-1}(0)$ determine $d-1$ pointwise independent sections of
$Hom(I_{\Psi_x}/I_{\Psi_x}^2, {\cal O}_{h^{-1}(0)})$. Thus
$I_{\Psi_x}/I_{\Psi_x}^2$ is free and ${\rm
Hom}(I_{\Psi_x}/I_{\Psi_x}^2, {\cal O}_{h^{-1}(0)})$, which is the
tangent space to  the Hilbert scheme (or Douady space) of
$h^{-1}(0)$, has dimension $d-1$. It follows that the Hilbert
scheme is smooth and of dimension $d-1$ at the point parametrizing
$\Psi_x$. By Definition 4.1 (ii), all subschemes $\Psi_x, x \in
h^{-1}(0)$ belong to one $T$-orbit in the Hilbert scheme. Thus
there exists a component $T'$ of the Hilbert scheme with an open
$T$-orbit $T'_o \subset T'$ such that $T'_o$ parametrizes
$\Psi_x$'s. But when $h^{-1}(0)_{\rm red}$ is not smooth, an orbit
of $T$ in $h^{-1}(0)$ is a torus by Definition 4.1 (i). This
implies that $T'_o= T'$ is a torus. Let $\hat{q}: {\cal U} \to T'$
and $p: {\cal U} \to h^{-1}(0)$ be the universal family morphisms
associated with the Hilbert scheme $T'$ such that fibers of
$\hat{q}$ are sent to $\Psi_x$'s by $p$. By Definition 4.1 (ii),
$p$ must be a bijective morphism. By Definition 4.1 (iii), $p$ is
unramified, hence it must be an isomorphism. Then $\hat{q}$
induces a $T$-equivariant morphism $q: h^{-1}(0) \to T'$ whose
fibers are $\Psi_x$'s. $\Box$

\medskip
Let $S$ be a normal surface and $j: S \to \Delta$ be a fibration
which is smooth over $\Delta \setminus \{0\}.$ Starting from $j: S
\to \Delta$, we can consider the following three operations to get
a new fibration $j': S' \to \Delta$.

(a) {\it Blow-up}: Assume that $S$ is smooth and $j^{-1}(0)_{\rm red}$
has a unique singular point $s$. Let $S'$ be the blow-up of $S$ at
$s$. The morphism $j':S' \to \Delta$ is just the composition of
the blow-up with $j.$

(b) {\it Contraction}: Assume that $S$ is smooth and there exists a
unique irreducible component $C$ of maximal multiplicity $m$  in
$j^{-1}(0)$ such that each connected component of $j^{-1}(0)
\setminus mC$ is  a Hirzebruch-Jung string (in the sense of [BHPV,
Section III.2]). Contract these Hirzebruch-Jung strings to get a
normal surface $S'$ with the induced morphism $j': S' \to \Delta$.

(c) {\it Cyclic cover}: For any positive integer $m$, we take the cyclic
cover $\nu: \Delta \to \Delta$ of degree $m$. Let $j': S'\to \Delta $ be the
normalization of the base change of $j$ by $\nu$.

\medskip
Note that the stable reduction of a given minimal elliptic
singular fibration is obtained by a finite number of operations of
the above three kinds. To generalize the stable reduction to
higher dimensions, we will explain how the above three operations
on surfaces can be generalized to higher dimensions.

\medskip
Let $j': S' \to \Delta$ be obtained from $j: S \to \Delta$ by one
of the  three operations (a), (b), (c) explained above. Given a
fibration $h: Z\to \Delta$  modeled on $j$ and a fixed choice of
an isomorphism $\iota: j^{-1}(0) \stackrel{\cong}{\to} \Psi_x
\subset h^{-1}(0)$ of Definition 4.1 (iii), we can construct a
fibration $h': Z' \to \Delta$ modeled on $j'$ in a canonical way
as follows.

(A) {\it Blow-up}:  If $j'$ is obtained by the operation (a), then we
blow-up $Z$ along the compact $T$-orbit $T \cdot (\iota(s))$
to get $h': Z' \to \Delta$, which is clearly a
fibration modeled on $j'$.

(B) {\it Contraction}: If $j'$ is obtained by the operation (b),
let $C_1, \cdots, C_k$ be the Hirzebruch-Jung strings in
$j^{-1}(0)$. Then the union of $T \cdot (\iota(C_1)), \cdots, T
\cdot (\iota(C_k))$  is a  divisor  in $Z$. We assign the
multiplicity on the component of this divisor by the multiplicity
of $C_i$ in $j^{-1}(0)$ and call the resulting divisor $D$. $D$
consists of some components of $h^{-1}(0)$. There exists a
fibration $f:D \to T'$ onto a complex torus $T'$ induced by  the
morphism in Proposition 4.2. Let $L$ be the line bundle on $Z$
corresponding to the Cartier divisor $D$ and $L^*$ be the dual
bundle of $L$. Since $C_i$'s are Hirzebruch-Jung strings, we have
(i) $L^*$ restricted to $D$ is $f$-ample and (ii) $R^1 f_*(L
^{*\otimes \ell}) = 0 $ for $\ell
>0$. Here (i) comes from the negative definiteness of the intersection
matrix of a Hirzebruch-Jung string and (ii) comes from the proof of
[BHPV, Proposition III (3.1)]. Thus this divisor $D$ can be
contracted to give a normal variety $Z'$ inducing $h': Z' \to
\Delta$ by [Fu, Theorem 2].  It is immediate to check that $h'$ is a
fibration modeled on $j'$.

 (C) {\it Cyclic cover}:  If $j'$ is obtained by the operation (c), then   $h': Z' \to \Delta$ is
defined to be the normalization of the base change of $h: Z \to
\Delta$ by $\nu: \Delta \to \Delta$. The $T$-action on $Z$
naturally lifts to a $T$-action on  $Z'$ and one can check that
$h'$ is a fibration modeled on $j'$.

\medskip
Since all the steps in the construction of the stable reduction
$j^{\sharp}$ from the minimal elliptic fibration $j$ with an
unstable fiber  are coming from the three operations (a), (b) and
(c), if we are given a fibration $h:Z \to \Delta$ modeled on $j$,
we get a fibration $h^{\sharp}: Z^{\sharp} \to \Delta$ modeled on
$j^{\sharp}$ by applying the three operations (A), (B) and (C).
The resulting fibration $h^{\sharp}$ is the stable reduction of
$h$. Now using the property (iii) of Definition 4.1 and the local
coordinate expression of the meromorphic correspondence for the
stable reduction for surfaces, we have the following.

\medskip
{\bf Proposition 4.3} {\it Let $h: Z \to \Delta$ be a fibration
modeled on a minimal elliptic fibration $j:S \to \Delta$ with an
unstable singular fiber. Let $h^{\sharp}: Z^{\sharp} \to \Delta$
be the stable reduction constructed above. The cyclic Galois group
$\langle \zeta_m \rangle$ of order $m$ acts on  $Z^{\sharp}$ with
fixed points. We can choose a coordinate system $(x, u_1, u_2,
\ldots, u_d)$ at a fixed point and a coordinate system $(y, v_1,
v_2, \ldots, v_d)$ on an open subset in $Z$ where $x$ (resp. $y$)
is the pull-back of the coordinate on $\Delta$ by $h^{\sharp}$
(resp. $h$) with the following meromorphic correspondence,
depending on the three cases of $j$.

(Case 1) $ y = x^2, v_1 = \frac{u_1}{x}, v_i=u_i$ for $2\leq i
\leq d$, and $$dy \wedge dv_1 \wedge \cdots \wedge d v_d = 2 dx
\wedge du_1 \wedge \cdots \wedge d u_d.$$

(Case 2) $ y = x^m, v_1 = \frac{u_1}{x^{m-1}}, v_i=u_i$ for $2\leq
i \leq d$, and $$dy \wedge dv_1 \wedge \cdots \wedge d v_d = m dx
\wedge du_1 \wedge \cdots \wedge d u_d.$$

(Case 3) $ y = x^m, v_1 = \frac{u_1}{x}, v_i=u_i$ for $2\leq i
\leq d$, and $$dy \wedge dv_1 \wedge \cdots \wedge d v_d = m
x^{m-2} dx \wedge du_1 \wedge \cdots \wedge d u_d.$$}

\section{Multiplicity in the unstable case}

\par
\noindent In this section, we will determine the possible values
of the multiplicity when the characteristic cycle is unstable. A
key  observation is the following proposition.

\medskip
{\bf Proposition  5.1} {\it Let $j: S \to \Delta$ be a fibration
of a smooth surface $S$. Let $h: Z \to \Delta$ be a fibration
modeled on $j: S\to \Delta$ in the sense of Definition 4.1 such
that $\omega_Z$ is trivial. Suppose $\chi: Z \cdots \to Z$ is a
bimeromorphic map on $Z$ inducing an isomorphism $h^{-1}(\Delta
\setminus \{0\}) \stackrel{\cong}{\to} h^{-1}(\Delta \setminus \{
0\})$ commuting with the projection $h$ and the $T$-action. Then
$\chi$ extends to a biholomorphic map $\tilde{\chi}: Z \to Z$.}

\medskip
{\it Proof}. Let $B \subset Z$ (possibly empty) be the
indeterminacy locus of the bimeromorphic map $\chi : Z \cdots \to
Z$. Let $\pi : \tilde{Z} \longrightarrow Z$ be a composition of
blow-ups such that the induced map $\tilde{\chi} :
\tilde{Z} \longrightarrow Z$ is holomorphic. As $\chi$ commutes
with $T$-action, it follows that each component of the exceptional
divisor $E = \cup_{i=1}^{m} E_i$ of $\pi$ is $T$-equivariant under
$\pi$. Since $T = \C^{d-1}, \dim Z = d+1$ and T acts on $Z$ with
discrete isotropy at every point, each $E_i$ is the blow-up of $T
\cdot P_i$ for some point $P_i$. On the other hand,
as $\omega_{Z} \cong {\cal O}_Z$
and $Z$ is smooth, it follows from [Ko, Lemma 4.3] that the bimeromorphic
map $\chi$
is isomorphic in codimension $1$.
Thus $\tilde{\chi}(E_i)$ is a subvariety of codimension
$\ge 2$. Moreover it is $T$-stable as $\tilde{\chi}$ is
$T$-equivariant. Thus $\tilde{\chi}(E_i) = T \cdot Q_i$ for some
point $Q_i$. As $\chi$ is $T$-equivariant, this means that any
fiber of $\pi$ is contracted by $\tilde{\chi}$. This implies that
$\chi$ itself is holomorphic. $\Box$

\medskip
 Now let $f:M \to B$ be a holomorphic Lagrangian fibration
over a complex manifold $B$.  Let $D \subset B$ be a component of
the discriminant of $f$ and $b \in D$ be a general point. Choose
an arc $\alpha: \Delta \to B$ with $\alpha(0) = b $ which
intersects $D$ transversally and denote by $h:Z \to \Delta$
 the pull-back of $f$ by $\alpha$.
 Then $Z$ is a complex manifold
 and the multiplicity of the fiber $h^{-1}(0)$ is the same as the
 multiplicity of $f^{-1}(b)$. The $T$-action on $M$ induces a $T$-action on $Z$ preserving
 the fibers of $h$.

\medskip
{\bf Proposition 5.2} {\it  Let $h: Z \to \Delta$ be the above
fibration    with a central fiber of multiplicity $n\geq 1$.
Assume that the characteristic cycle is unstable, i.e., belongs to
one of the three cases considered in Section 4. Let
$$\begin{array}{ccc} Z_1 & \stackrel{\pi_1}{\longrightarrow} & Z \\
h_1 \downarrow & & \downarrow h \\ \Delta &
\stackrel{\nu_1}{\longrightarrow}  &\Delta, \end{array} $$ be the
normalization of the fiber product where $\nu_1$ is a cyclic
covering of degree $n$, as in Proposition 3.1. By Proposition 2.3
and Proposition 2.5, $h_1$ is a fibration modeled on a minimal
elliptic fibration with an unstable fiber, in the sense of
Definition 4.1. Let $h_2 : Z_2 \to \Delta$ be the stable reduction
of $h_1: Z_1 \to \Delta$ constructed in Section 4 with a dominant
meromorphic map $\pi_2: Z_2 \cdots\to Z_1$ of degree $m$. Let
$\Gamma_{12}$ be the cyclic group of order $mn$ with a subgroup
$\Gamma_2 \subset \Gamma_{12}$ of order $m$ and the quotient group
$\Gamma_1= \Gamma_{12}/\Gamma_2$ of order $n$. Then there exists
an action of $\Gamma_{12}$ on $Z_2$ commuting with the $T$-action
and compatible with the fibration $h_2$, such that the induced
action of $\Gamma_2$ on $Z_2$ agrees with the Galois action of the
cyclic group of order $m$ on the stable reduction. Via $\pi_2$,
the action of $\Gamma_{12}$ on $Z_2$ induces an action of
$\Gamma_1$ on $Z_1$, which agrees with the unramified Galois
action of the cyclic group of order $n$ on $Z_1$. In particular,
the induced $\Gamma_1$-action on $\omega_{Z_1}$ is trivial. }

\medskip
{\it Proof}. We have a cyclic covering of degree $m$, $\nu_2:
\Delta \to \Delta$ and the commuting diagram
$$\begin{array}{ccccc} Z_2 & \stackrel{\pi_2}{\cdots\rightarrow} & Z_1 &
\stackrel{\pi_1}{\longrightarrow} & Z \\ h_2 \downarrow & & h_1
\downarrow & & \downarrow h \\ \Delta &
\stackrel{\nu_2}{\longrightarrow} & \Delta &
\stackrel{\nu_1}{\longrightarrow} & \Delta. \end{array} $$
  Let $\hat{h}: \hat{Z} \to \Delta$ be the
normalization of the base change of $h: Z \to \Delta$ by the
cyclic covering $\hat{\nu}: \Delta \to \Delta$ of degree $mn$. By
the construction of the stable reduction, there exists a
bimeromorphic map $\varphi: \hat{Z} \cdots \to Z_2$ commuting with
the fibrations $\hat{h}$ and $h_2$, inducing a biholomorphism
outside the central fibers. We have the Galois action of the
cyclic group $\Gamma_{12}$ on  $\hat{Z}$ respecting the fibration
$\hat{h}$. Thus we have a bimeromorphic action of $\Gamma_{12}$ on
$Z_2$ respecting $h_2$ which induces a biholomorphic action on
$Z_2 \setminus h_2^{-1}(0)$. By Proposition 5.1, this extends to a
biholomorphic action of $\Gamma_{12}$ on $Z_2$. By the
construction,  a cyclic subgroup $\Gamma_2 \subset \Gamma_{12}$
acts on $Z_2$ as the Galois action for the cyclic covering of
degree $m$ in the construction of the stable reduction of $h_1:
Z_1 \to \Delta$.  The action of the quotient group $\Gamma_1$ on
$Z_1$ must be the Galois action on $Z_1$ induced by $\nu_1$, which
preserves $\omega_{Z_1}$ from Proposition 3.1. $\Box$

\medskip
Now we determine the possible values of $n$. We will treat the
three cases of Section 5 separately.

\medskip
{\bf Proposition 5.3} {\it In the setting of Proposition 5.2,
suppose that the characteristic cycle  is of type $I^*_b, b \geq 1$.
Then the  multiplicity $n$ must be 1.}

\medskip
{\it Proof}. By Proposition 5.2, we have an action of the cyclic
group $\Gamma_{12}$ generated by the root of unity $ \zeta_{mn}$
on $Z_2$ with $m=2$. We claim that $\zeta_{mn}$ acts on the
dualizing sheaf of the fiber $\omega_{F_2}$ by $-1$. By the
construction of the stable reduction, the generator $\zeta_m =
\zeta_{mn}^n$ of the Galois group $\Gamma_2$ of the covering
$\nu_2$ acts on the dualizing sheaf $\omega_{Z_2}$ by $-1$. Thus
${\rm Aut}_T(F_2)/{\rm Aut}^{\omega}_T(F_2)$ is non-trivial and
by Proposition 3.2,
$${\rm Aut}_T(F_2)/{\rm Aut}^{\omega}_T(F_2) \cong {\bf Z}/2.$$ Let
$\bar{\zeta}_{mn} \in {\rm Aut}_T(F_2)/{\rm Aut}^{\omega}_T(F_2)$
be the image of $\zeta_{mn} \in {\rm Aut}_T(F_2)$ in the quotient
group. Since $\bar{\zeta}_{mn}^n = -1,$ $\bar{\zeta}_{mn}$ must be
$-1$, too. This means that it acts on $\omega_{F_2}$ by $-1$.
Since its weight on $\omega_{\Delta}$ is $\zeta_{mn}$, it acts on
$\omega_{Z_2}$ by $- \zeta_{mn}$. By Proposition 4.3, the induced
action of $\zeta_{mn}$ on $\omega_{Z_1}$ is also by $-
\zeta_{mn}$. But by Proposition 5.2, the induced action on
$\omega_{Z_1}$ must be trivial. Hence $-\zeta_{mn} =1$, which
implies that $n=1$. $\Box$

\medskip
{\bf Proposition 5.4} {\it In the setting of Proposition 5.2,
suppose that the characteristic cycle  is of type $I^*_0, II^*,
III^*$ or $ IV^*$. Then the possible values of $(m,n)$ are $(2,2),
(2,3), (3,2)$. }

\medskip {\it Proof}. By Proposition 5.2, we have an action of the
cyclic group $\Gamma_{12}$ generated by the root of unity $
\zeta_{mn}$ on $Z_2$. The generator $\zeta_{mn}$ acts on the
dualizing sheaf of the fiber $\omega_{F_2}$ by some weight
$\zeta_{mn}^a$ where $0 \leq a < mn$ is some integer.
 By the
construction of the stable reduction,  the generator $\zeta_m =
\zeta_{mn}^n$ of the Galois group $\Gamma_2$ of the covering
$\nu_2$ acts on the dualizing sheaf $\omega_{Z_2}$ by
$\zeta_m^{m-1}$. It follows that $a \equiv m-1 \; \mod \; m.$

Since $\zeta_{mn}$ acts on $\omega_{\Delta}$ by $\zeta_{mn}$, the
action on $\omega_{Z_2}$ is by $\zeta_{mn}^{a+1}$. By Proposition
4.3, the induced action on $\omega_{Z_1}$ is also by
$\zeta_{mn}^{a+1}$. But by Proposition 5.2, the induced action on
$\omega_{Z_1}$ must be trivial. It follows that $a+1 \equiv 0 \;
\mod \; mn$ and $a = mn-1$. Thus the action on $\omega_{F_2}$ is
by $\zeta_{mn}^{mn-1}$. By Proposition 3.2, $\zeta_{mn}^{mn-1} =
\zeta_{mn}^{-1}$ must have order $2,3,4,$ or $6.$ The possible
values of $(m,n)$ with $m, n \geq 2$ are $(2,2), (2,3), (3,2)$.
$\Box$

\medskip
{\bf Proposition 5.5} {\it In the setting of Proposition 5.2,
suppose that the characteristic cycle  is of type $II, III$, or $
IV$.   Then the possible values of $(m,n)$ are $(3,2), (3,4), (4,3)$
and $(6,5)$. }

\medskip
{\it Proof}. By Proposition 5.2, we have an action of the cyclic
group $\Gamma_{12}$ generated by the root of unity $ \zeta_{mn}$
on $Z_2$. The generator $\zeta_{mn}$ acts on the dualizing sheaf
of the fiber $\omega_{F_2}$ by some weight $\zeta_{mn}^a$ where $0
\leq a < mn$ is an integer.
 By the
construction of the stable reduction, the generator $\zeta_m =
\zeta_{mn}^n$ of the Galois group of the covering $\nu_2$  acts on
the dualizing sheaf $\omega_{Z_2}$ by $\zeta_m$. It follows that
$a \equiv 1 \; \mod \; m.$

Since $\zeta_{mn}$ acts on $\omega_{\Delta}$ by $\zeta_{mn}$, the
action on $\omega_{Z_2}$ is by $\zeta_{mn}^{a+1}$. By Proposition
4.3, the induced action on $\omega_{Z_1}$ is  by $\zeta_{mn}^{m-2}
\cdot \zeta_{mn}^{a+1} = \zeta_{mn}^{a + m-1}$. But by Proposition
5.2, the induced action on $\omega_{Z_1}$ must be trivial. It
follows that $a +m-1 \equiv 0 \; \mod \; mn.$ From  $a \equiv 1 \;
\mod \; m,$ let $a = mq +1$. Then $q+1 \equiv 0 \; \mod \; n.$
Setting $q+1 = n q',$
$$ a = m(nq'-1) + 1 \equiv 1-m \; \mod \; mn.$$ It follows that
the action on $\omega_{F_2}$ is by $\zeta_{mn}^{1-m}$. By
Proposition 3.2, $\zeta_{mn}^{1-m}$ must have order $2,3,4,$ or
$6.$ The possible values of $(m,n)$ with $m=3,4,6$ and $ n \geq 2$
are $(3,2), (3,4), (4,3)$ or $(6,5)$. $\Box$

\medskip
Proposition 5.3, Proposition 5.4, Proposition 5.5 and Proposition
3.5 complete the proof of Theorem 1.1, modulo the realization
part, which will be given in the next section.

\section{Explicit examples of multiple fibers}
\par
\noindent In this section, we show the last statement of Theorem
1.1 by giving explicit examples of multiple fibers described in
(i) -(v) there. It suffices to construct them in dimension $4$.

\medskip

We shall first give an explicit example of (ii) in Theorem 1.1
with full details.
In fact, this is one of the cases missing in [Ma]
and the discovery of this case is one of the starting points of
our Theorem 1.1.

\medskip
{\bf Example 6.1 (Type II with multiplicity $5$)}

Let $\varphi : S \to \Delta_t$ be a relatively minimal
 elliptic surface over the unit disk $\Delta_t$ given by
the Weierstrass equation
$$y^2 = x^3 + t\, .$$
There is a very convenient algorithm, called the Tate's algorithm,
to determine the type of
singular fiber from the Weierstrass equation. This algorithm
is given by [Ta, Pages 34-35, Summary 0 with additional definition
in Page 36 (3.6)]. Applying to our equation,
we readily see that $\varphi$ has a singular fiber of Type $II$
over $t=0$. (This algorithm can be also used to determine
the singular fiber in Example 6.2 below.) We also note that
the $2$-form
$$\sigma_S := \frac{dx \wedge dt}{y}$$
gives the generator of $\omega_S \simeq {\mathcal O}_S$. Note also
that $S$
admits an automorphism given by
$$\tau^{*}(x, y, t) = (\zeta_5^2x, \zeta_5^3y, \zeta_5t)\, .$$
Now consider the product $4$-fold
$$M := S \times E \times \Delta_s$$
where $E$ is an elliptic curve and $\Delta_s$ is a unit disk. $M$
has a fibration
$$f : M \to \Delta_t \times \Delta_s\, ;\,
((x,y,t), z, s) \mapsto (t, s)\, .$$ The $2$-form
$$\sigma_M := \sigma_S + dz \wedge ds$$
is a symplectic form on $M$ and makes $f$ Lagrangian. We define
the automorphism $g$ of $M$ by
$$g^*((x,y,t), z, s) = ((\zeta_5^2x, \zeta_5^3y, \zeta_5t), z + p, s)$$
where $p$ is a $5$-torsion point on $E$. Then, $\langle g \rangle
\simeq {\mathbf Z}/5$ and $\langle g \rangle$ acts on $f : M
\longrightarrow \Delta_t \times \Delta_s$ freely. Moreover
$g^{*}\sigma_{M} = \sigma_{M}$ by the explicit form of $g$. The
quotient manifold (where $u = t^5$)
$$M/\langle g \rangle \longrightarrow \Delta_{u} \times \Delta_s$$
is then a Lagrangian fibration with multiple fibers of
multiplicity $5$ along $u = 0$ whose characteristic cycles are of
Type $II$. Note that the fiber $F = f^{-1}(0,s) \simeq
f^{-1}(0) \times E$ is stable under $g$ and satisfies
$$g^* (\frac{dx}{y} \wedge dz) = \zeta_5^{-1} \frac{dx}{y} \wedge dz\, ,$$
for the generator of $\omega_F$. We also note that the stable reduction
of $S$ is $y^2 = x^3 + 1$. Geometrically, the stable reduction is the second projection from the product $E_{\zeta_3} \times \Delta$
of the elliptic curve $E_{\zeta_3}$ and the unit disk.

\medskip

{\bf Example 6.2 (Examples of (i), (iii), (iv), (v) when characteristic cycle
is not of Type $I_{2m}$ ($1 \le m \le \infty$))}

The construction for other cases (i), (iii), (iv), (v) in Theorem
1.1 are quite similar if the characteristic cycle is not of Type
$I_{m}$. In fact, in the quotient $(S \times E \times
\Delta_s)/\langle g \rangle$ in the example above, we just replace
the pair $(S, \tau, g)$ as follows (with the same $E$ and
$\Delta_s$ and the same expression for $\sigma_S$ and $\sigma_M$),
according to the cases. Here $S$ is given by the Weierstrass
equation and $p_n$ is an $n$-torsion point of the elliptic curve
$E$:

(i) ($n =6$ and $E_{\zeta_3}$):
$$S : y^2 = x^3 +1\,\, ,\,\, \tau^*(x, y, t) =
(\zeta_6^2 x, \zeta_6^3 y, \zeta_6 t)$$
$$g^*((x,y, t), z, s) =
((\zeta_6^2 x, \zeta_6^3 y, \zeta_6 t), z + p_6, s)\, .$$

(iii) ($n =4$ and $E_{\zeta_4}$) :
$$S : y^2 = x^3 + x\,\, ,\,\,\tau^*(x, y, t) =
(-x,\zeta_4^3 y, \zeta_4 t)$$
$$g^*((x,y, t), z, s) = ((-x,\zeta_4^3 y, \zeta_4 t), z + p_4, s)\, .$$

(iii) ($n =4$ and Type $IV$) :
$$S : y^2 = x^3 + t^2\,\, ,\,\, \tau^*(x, y, t) = (-x,
\zeta_4^3 y, \zeta_4 t)$$
$$g^*((x,y, t), z, s) = ((-x,\zeta_4^3 y, \zeta_4 t), z + p_4, s)\, .$$

(iv) ($n =3$ and $E_{\zeta_3}$):
$$S : y^2 = x^3 +1\,\, ,\,\, \tau^*(x, y, t) =
(\zeta_3^2 x, y, \zeta_3 t)$$
$$g^*((x,y, t), z, s) = ((\zeta_3^2 x, y, \zeta_3 t), z + p_3, s)\, .$$

(iv) ($n =3$ and Type $III$):
$$S : y^2 = x^3 +tx\,\, ,\,\, \tau^*(x, y, t) =
(\zeta_3^2 x, y, \zeta_3 t)$$
$$g^*((x,y, t), z, s) = ((\zeta_3^2 x, y, \zeta_3 t), z + p_3, s)\, .$$

(iv) ($n =3$ and Type $I_0^*$):
$$S : y^2 = x^3 +t^3\,\, ,\,\, \tau^*(x, y, t) =
(\zeta_3^2 x, y, \zeta_3 t)$$
$$g^*((x,y, t), z, s) = ((\zeta_3^2 x, y, \zeta_3 t), z + p_3, s)\, .$$

(v) ($n =2$ and Type $I_0$) :
$$S : y^2 = x^3 + ax + b\,\, ,\,\, \tau^*(x, y, t) = (x,
-y, -t)$$
$$g^*((x,y, t), z, s) = ((x,-y, -t), z + p_2, s)\, .$$
Here $a$ and $b$ are any complex numbers such that $4a^3 + 27b^2 \not= 0$.

(v) ($n =2$ and Type $I_0^*$):
$$S : y^2 = x^3 + t^2x\,\, ,\,\, \tau^*(x, y, t) = (x,-y, -t)$$
$$g^*((x,y, t), z, s) = ((x,-y, -t), z + p_2, s)\, .$$

(v) ($n =2$ and Type $IV$):
$$S : y^2 = x^3 + t^2\,\, ,\,\, \tau^*(x, y, t) = (x, -y,
-t)$$
$$g^*((x,y, t), z, s) = ((x,-y, -t), z + p_2, s)\, .$$

(v) ($n =2$ and Type $IV^{*}$):
$$S : y^2 = x^3 + t^4\,\, ,\,\, \tau^*(x, y, t) = (x,
-y, -t)$$
$$g^*((x,y, t), z, s) = ((x,-y, -t), z + p_2, s)\, .$$

\medskip

{\bf Example 6.3 (Examples of $n=2$ and Type $I_{2m}$ ($1 \le m \le \infty$))}

We will use the example in [HO, Proposition 4.13]. Let us recall
the setting.  Let $R_{k} = {\rm Specan}\, {\mathbf
C}[u^{k+1}v^{-1}, u^{-k}v]$ ($k \in {\mathbf Z}$).  There is a
natural morphism $g_{k} : R_{k} \longrightarrow {\rm Specan}\,
{\mathbf C}[u]$. Let $E$ be an elliptic curve. Using the morphisms
$g_{k}$, which are compatible with the natural gluing of the
spaces $R_{k}$, we obtain a morphism
$$(\cup_{k \in {\mathbf Z}} R_{k}) \times E \times {\rm Specan}\, {\mathbf C}[y] \longrightarrow {\rm Specan}\, {\mathbf C}[u, y]\,\, .$$
Restricting this morphism over a sufficiently small
$2$-dimensional disk $\Delta^2$ (centered at $(u, y) = (0,0)$), we
obtain a fibration
$$\tilde{f} : \tilde{M} \longrightarrow \Delta^{2}_{(u, y)}\,\, .$$
The fiber over $u = 0$ is an infinite  chain of ${\mathbf P}_1
\times E$, while the fiber over $u \not= 0$ is ${\mathbf C}^{*}
\times E$. Let $\alpha$ be a  point of $E$. Then ${\mathbf Z}$
acts on $\tilde{f} : \tilde{M} \longrightarrow \Delta^2$ if we
define the action of $m \in {\mathbf Z}$ by
$$(u^{k+1}v^{-1}, u^{-k}v, x, y) \mapsto (u^{k+1+m}v^{-1}, u^{-k-m}v, x + m\alpha, y)\,\, .$$
As explained in [HO], this action is properly discontinuous and
free. Let $p$ be a $2$-torsion point of $E$. Define the
automorphism of $\tilde{M}$ by
$$g^*(u, v, x, y) = (-u, v^{-1}, x+p, y)\, .$$
$g$ does {\it not} commute with the ${\mathbf Z}$-action, but it
commutes with the action of the index two subgroup $2{\mathbf Z}$.
This is why our example below of fibers of multiplicity 2 can be
constructed for Type $I_b$ with even $b$ even, but not with odd
$b$.

First choose $\alpha$ to be a nontorsion point of $E$ and set $M_2
= \tilde{M}/2{\mathbf Z}$.  The symplectic $2$-form
$$du \wedge \frac{dv}{v} + dx \wedge dy$$
on $\tilde{M}$ descends to a symplectic $2$-form on $M_2$.
We regard $M_2$ as a symplectic manifold with this symplectic from.
Then, $g$ descends
to the free symplectic involution of the Lagrngian fibration
$$f : M_2 \longrightarrow \Delta^2\, .$$
Thus, the quotient fibration
$$\overline{f} : M_2/\langle g \rangle \longrightarrow \Delta^2$$
gives an example of a Lagrangian fibration such that a general singular
fiber is of multiplicity
$2$ with characteristic cycle of Type $A_{\infty} = I_{\infty}$.

In the above construction, if we choose $\alpha \in E$ to be a
torsion element of order $1 \le 2\ell < \infty$ and choose the
$2$-torsion point $p$ with $p \not\in \langle \alpha \rangle$, we
obtain an example of a Lagrangian fibration such that a general
singular fiber is of multiplicity $2$ with characteristic cycle of
Type $I_{2\ell}$ ($\ell \ge 1$).

\medskip

(6.1)-(6.3) complete the proof of the realizability.

\bigskip
{\bf References}

\medskip

[BHPV] Barth, W., Hulek, K., Peters, C., Van de Ven, : {\it
Compact complex surfaces}. Second enlarged edition. Springer
Verlag, Berlin-Heidelberg, 2004

[Fu] Fujiki, A.: On the blowing down of analytic spaces.
Publ. RIMS, Kyoto Univ. {\bf 10} (1975)  473-507

[HO] Hwang, J.-M., Oguiso, K.: Characteristic foliation on the
discriminant hypersurface of a holomorphic Lagrangian fibration.
 Amer. J. Math. {\bf 131} (2009) 981-1007

[Kd] Kodaira, K.: On compact analytic surfaces. II, Ann. of Math.
(2) {\bf 77} (1963) 563--626.

[Ko] Koll\'ar, J.: Flops. Nagoya Math. J. {\bf 113} (1989) 15--36.

[Ma] Matsushita, D.: A canonical bundle formula for projective
Lagrangian fibrations. preprint, 2007, arXiv:0710.0122.

[Mc] McMullen C. T.: Dynamics on blowups of the projective plane,
Publ. Math. Inst. Hautes \'Etudes Sci. {\bf 105} (2007) 49--89.

[Ta] Tate, J.: Algorithm for determining the type of a singular fiber in an elliptic pencil. Modular functions of one variable, IV, Lecture
Notes in Math. {\bf 476} (1975) 33--52.

\pagebreak

Jun-Muk Hwang

Korea Institute for Advanced Study

Hoegiro 87, Seoul 130-722, Korea

 jmhwang@kias.re.kr

\bigskip

Keiji Oguiso

Department of Mathematics, Osaka University

Toyonaka 560-0043 Osaka, Japan

 oguiso@math.sci.osaka-u.ac.jp

\end{document}